\documentclass[12pt]{article}
\usepackage{hyperref}
\usepackage{amsfonts}
\usepackage{amsmath}
\usepackage{amssymb}
\newtheorem{theorem}{Theorem}

\newtheorem{lemma}[theorem]{Lemma}
\newtheorem{corollary}[theorem]{Corollary}

\newtheorem{remark}[theorem]{Remark}
\newcommand{\R}{\mathbb{R}}

\newcommand{\hess}{\operatorname{hess}}
\newcommand{\Hess}{\operatorname{Hess}}

\newcommand{\po}{{\hspace*{-1ex}}{\bf .  }}

\newcommand{\enne}{\mathbb{N}}
\newcommand{\g}[2]{\mbox{$\langle #1 ,#2 \rangle$}}
\newcommand{\rf}[1]{\mbox{(\ref{#1})}}
\newcommand{\rl}[1]{{~\ref{#1}}}
\def\<{\langle}

\def\>{\rangle}
\def\bea{\begin{eqnarray*} }
\def\eea{\end{eqnarray*} }
\def\be{\begin{equation} }
\def\ee{\end{equation} }

\def\nap{\bar\nabla}
\def\proof{\noindent{\it Proof: }}
\newcommand{\pair}[1]{\langle#1\rangle}

\def\qed{\ifhmode\unskip\nobreak\fi\ifmmode\ifinner
\else\hskip5 pt \fi\fi\hbox{\hskip5 pt \vrule width4 pt
height6 pt  depth1.5 pt  \hskip 1pt }}
\begin{document}

\title{Higher order mean curvature estimates for bounded complete hypersurfaces
\footnote{This work was partially supported by MINECO (Ministerio de Econom\'{\i}a
y Competitividad) and FEDER (Fondo Europeo de Desarrollo Regional)
project MTM2012-34037 and Fundaci\'{o}n S\'{e}neca project
04540/GERM/06, Spain. This research is a result of the
activity developed within the framework of the Program in Support
of Excellence Groups of the Region de Murcia, Spain, by
Fundamental S\'{e}neca, Regional Agency for Science and Technology
(Regional Plan for Science and Technology 2007-2010).
M. Dajczer was partially supported by CNPq and FAPERJ, Brazil, and MEC Grant
SAB2011-0152.}}

\author {L. J. Al\'{i}as, M. Dajczer and M. Rigoli}
\date{}
\maketitle

\begin{abstract} We obtain sharp estimates involving the mean curvatures of
higher order of a complete bounded hypersurface immersed in a complete Riemannian
manifold. Similar results are also given for complete spacelike hypersurfaces
in Lorentzian ambient spaces.
\end{abstract}

Estimates for the $k$-mean curvatures $H_k$  of higher order of a compact
hypersurface in a complete Riemannian manifold have been subsequently
obtained~by \mbox{Vlachos} \cite{Vl}, Veeravalli \cite{V}, Fontenele-Silva
\cite{FS1},  Roth \cite{R} and Ranjbar-Motlagh \cite{Rm}. In this paper,
we generalize a result given in the latter that we describe next.

Let $f\colon M^n\to\bar{M}^{n+1}$ be a codimension one isometric immersion
between complete Riemannian manifolds.  Assume that the hypersurface lies inside
a closed geodesic ball $B_{\bar{M}}(r)$ of radius $r$ and center $o\in\bar{M}^{n+1}$
and that \mbox{$0<r<\min\{{\rm inj}_{\bar{M}}(o), \pi/2\sqrt{b}\}$} where ${\rm
inj}_{\bar{M}}(o)$ is the injectivity radius at $o$ and $\pi/2\sqrt{b}$ is
replaced by $+\infty$ if $b\leq 0$.
Suppose also that there is a point $p_0\in M^n$ such that
$f(p_0)\in S_{\bar{M}}(r)$ where $S_{\bar{M}}(r)$ is the boundary of
$B_{\bar{M}}(r)$.
In the context of this paper, this is a slightly weaker assumption than asking
$M^n$ to be compact.
Let $K_{\bar{M}}^{\mathrm{rad}}$ denote the radial sectional curvatures in
$B_{\bar{M}}(r)$ along geodesics issuing from the center and assume that
$K_{\bar{M}}^{\mathrm{rad}}\leq b$ for some constant $b\in\R$. Assume also that
$H_{k+1}\neq 0$ everywhere  for some $2\leq k\leq n-1$. In this situation,
it turns that the $p_0$ is an elliptic point. This means that the second
fundamental form of $f$ at $p_0$ with respect to the inner pointing orientation
is positive definite.  From the well-known Garding
inequalities it follows that $H_j>0$ for $1\le j\le k+1$.

In the above situation, it was shown in Theorem 4.2 in \cite{Rm} that
$$
\sup_M\left(\frac{H_{j+1}}{H_j}\right)\geq C_b(r)
$$
for any $1\le j\leq k$, where the constant $C_b(r)$  given by (\ref{cbr})
below is the mean curvature of a geodesic  sphere of radius $r$ in a simply
connected space form  of sectional curvature $b$.  Moreover, if  equality holds
for some $j$ then it follows that $M^n=S_{\bar{M}}(r)$.

Our main goal in this paper is to replace the assumption of compactness
of the submanifold by the much weaker of completeness. The tool that
makes this generalization possible is an Omori-Yau type maximum  principle  for
trace type differential operators in the spirit of those due to Albanese, Al\'{i}as and Rigoli
\cite{AAR} (see Theorem\rl{MAIN} below).

The following is a consequence of the quite more general result given in
Section~$2$. Here, the more general but technical assumptions made in
Theorem \ref{main} of Section $2$ take a simpler geometric form.

\begin{theorem}\po\label{main3} \!\!Let $f\colon M^n\to\bar{M}^{n+1}$ be an isometric
immersion between complete Riemannian manifolds such that $f(M)\subset
B_{\bar{M}}(r)$. Assume $H_{k+1}\neq 0$ everywhere for some $2\leq k\leq n-1$ and that
the sectional curvatures satisfy  $K_M\geq K>-\infty$
and $K_{\bar{M}}^{\mathrm{rad}}\leq b$
for some constant $b\in\R$.  If $f$ has an elliptic point,
then
\be\label{ineq3}
\sup_M \sqrt[j+1]{H_{j+1}}\geq\sup_M\left(\frac{H_{j+1}}{H_j}\right)\geq
C_b(r),\quad 1\le j\leq k.
\ee
Moreover, if there exists a point  $p_0\in M^n$ such that
$f(p_0)\in S_{\bar{M}}(r)$ and $\sup_M\left({H_{j+1}}/{H_j}\right)=C_b(r)$ for
some $j$   then $M^n=S_{\bar{M}}(r)$.
\end{theorem}

In the second part of the paper and motivated,  among others, by the
results in \cite{AHP} and \cite{AIRlorentz}, we
show that similar estimates than in the Riemannian case
hold for complete spacelike hypersurfaces in
Lorentzian ambient spaces.

\section{A maximum principle}
The aim of this section is to introduce the main analytic ingredient
for the proof of our results. It consists in a maximum  principle of
Omori-Yau type in the spirit of those given in \cite{AAR} that applies to trace type
operators like those described in the sequel.\vspace{1,5ex}

Let $M^n$ be a Riemannian manifold and $\nabla$  the Levi-Civita
connection. For  $u\in C^{2}(M)$ let $\hess{u}\colon
TM\rightarrow TM$ be the symmetric operator given by
$\hess{u}(X)=\nabla_X\nabla u$ and by $\Hess{u}\colon TM\times
TM\rightarrow C^{0}(M)$ the metrically equivalent bilinear form
given by
$$
\Hess{u}(X,Y)=\pair{\hess{u}(X),Y}.
$$
Associated to a symmetric tensor $P\colon TM\rightarrow TM$, we
consider the second order differential operator
$L\colon\mathcal{C}^{2}(M)\rightarrow\mathcal{C}^{0}(M)$ given by
$L=\mathrm{Tr}\,(P \circ \hess)$. Observe that
$
L(u)=\mathrm{div}(P\nabla u)-\<\mathrm{div}P,\nabla u\>,
$
where $\mathrm{div}P=\mathrm{Tr}\nabla P$. This implies that $L$
is (semi-)elliptic if and only if $P$ is
positive (semi-)definite. The following result is Theorem B together with
Remark 1.2 in \cite{AAR}.

\begin{theorem}\label{maxprinc}\po
Let $M^n$ be a Riemannian manifold and let $L=\mathrm{Tr}\,(P \circ
\hess)$ be a~semi-elliptic linear operator. Let
$q\in\mathcal{C}^0(M)$ be nonnegative such that $q>0$ outside a compact
set. Assume that there exists
$\gamma\in C^2(M)$ with the following properties:
\begin{itemize}
\item[(a)] $\gamma(p)\rightarrow +\infty$ as $p\rightarrow \infty$,
\item[(b)] $\|\nabla\gamma\|\leq G(\gamma)$ off a compact set,
\item[(c)] $qL\gamma\leq G(\gamma)$ off a compact set
\end{itemize}
where $G$ is a smooth function on $[0,+\infty)$ such that:
$$
\mathrm{(i)}\  G(0)>0, \quad \mathrm{(ii)}\  G'(t)\geq 0\;\; \mbox{and}\;\;
\mathrm{(iii)}\  1/G(t)\not \in L^1(+\infty).
$$
Then, for any function $u \in C^{2}(M)$ with $u^*=\sup_{M}u < +\infty$
there exists a sequence $\left\{p_{j}\right\}_{j\in \enne}$ in $M^n$
such that
$$
\mathrm{(i)}\ u(p_j)>u^*-\frac{1}{j},\;\;
\mathrm{(ii)} \ \|\nabla u(p_j)\|<\frac{1}{j}\;\;\mbox{and}\;\;
\mathrm{(iii)} \ q(p_j)Lu(p_j)< \frac{1}{j}.
$$
\end{theorem}

Following the terminology in \cite{AAR}, we say that
the $q$-\emph{Omori-Yau maximum principle holds on $M^n$ for $L$} as
above whenever the conclusions of Theorem~\ref{maxprinc}
hold.\vspace{1ex}

Let $M^n$ be a complete noncompact Riemannian manifold. Denote by
$r(x)$ the distance function to a fixed reference point $o\in M^n$.
Then $r(x)$ satisfies assumptions $(a)$
and (b) of Theorem \ref{maxprinc}. Although $r(x)$ is not $C^2$ in $o$
and its cut locus $\mathrm{cut}(o)$, one could think of it as a natural candidate
for $\gamma$, under appropriate curvature assumptions.
The technical difficulty arising from this choice, and related to
the lack of smoothness, forces us to
introduce a reasoning in some way similar to approaching the problem
via viscosity solutions in order to get the following result.

\begin{theorem}\po\label{MAIN}
Let $M^n$ be a complete, non-compact Riemannian manifold and let
$r(x)$ be the Riemannian distance function from a reference point $o\in M^n$.
Assume that the sectional curvature of $M^n$ satisfies
\be\label{A1}
K_M(x)\geq -G^2(r(x))
\ee
with $G\in\mathcal{C}^1([0,+\infty))$ satisfying
$$
\mathrm{(i)}\  G(0)>0, \quad \mathrm{(ii)}\  G'(t)\geq 0\;\;\; \mbox{and}\;\;\;
\mathrm{(iii)}\  1/G(t)\not \in L^1(+\infty).
$$
Then, the \mbox{$q$-Omori-Yau} maximum principle holds on $M^n$ for any
semi-elliptic operator of the form $L=\mathrm{Tr}\,(P\circ\hess)$
with $\mathrm{tr}P>0$ on $M^n$ where $q=1/\mathrm{Tr}P$.
\end{theorem}

\proof Let $D_o=M^n\setminus\mathrm{cut}(o)$ be the domain of
normal geodesic coordinates centered at $o$.  On $D_o$ we have
from \rf{A1} and the general Hessian comparison theorem \cite[Theorem 2.3]{PRSbook}
 that
\be\label{A3}
\Hess(r)\leq\frac{g'(r)}{g(r)}(\g{}{}-dr\otimes dr),
\ee
where $g(t)$ is the (positive on $\mathbb{R}^+=(0,+\infty)$) solution
of the Cauchy problem
\be\label{A4}
\left\{\begin{array}{l}
g''(t)-G^2(t)g(t)=0,\\
g(0)=0, \quad g'(0)=1.
\end{array}\right.
\ee
Letting
$$
\psi(t)=\frac{1}{G(0)}\big(e^{\int_0^tG(s)ds}-1\big)
$$
we have $\psi(0)=0$, $\psi'(0)=1$ and
$$
\psi''(t)-G^2(t)\psi(t)=\frac{1}{G(0)}\left(G^2(t)+G'(t)\,e^{\int_0^tG(s)ds}\right)\geq 0,
$$
that is, $\psi$ is a subsolution of \rf{A4}. By the Sturm comparison theorem
$$
\frac{g'(t)}{g(t)}\leq\frac{\psi'(t)}{\psi(t)}
=G(t)\frac{e^{\int_0^{t}G(s)ds}}{e^{\int_0^{t}G(s)ds}-1}.
$$
Thus, we have
$$
Lr(x)\leq \mathrm{tr}P(x)\frac{\psi'(r(x))}{\psi(r(x))}
=\mathrm{tr}P(x)G(r(x))\frac{e^{\int_0^{r(x)}G(s)ds}}{e^{\int_0^{r(x)}G(s)ds}-1}.
$$
Since $G>0$, $G'\geq 0$ and $\mathrm{tr}P\geq 0$, we obtain
\be\label{A8}
Lr(x)\leq \mathrm{tr}P(x)G(r(x)+1)\frac{e^{\int_0^{r(x)}G(s)ds}}{e^{\int_1^{r(x)}G(s)ds}-1},
\;\; r(x)\geq 2.
\ee
Define
\be\label{A10}
\varphi(t)=\int_0^{t}\frac{ds}{G(s+1)}
\ee
so that
$$
\varphi'(t)=\frac{1}{G(t+1)} \;\;\text{and}\;\; \varphi''(t)\leq 0.
$$
Set $\gamma(x)=\varphi(r(x))$ on $M^n\setminus \bar{B}_2$ and
note that
\be\label{A15}
\gamma(x)\rightarrow+\infty \text{ as } x\rightarrow\infty
\ee
because
$\varphi(t)\rightarrow+\infty$ as $t\rightarrow+\infty$ since $1/G\notin L^1(+\infty)$.

Using the formula
$
L\varphi(u)=\varphi'(u)Lu+\varphi''(u)\g{P\nabla u}{\nabla u}
$
and that $P$ is positive semi-definite,  we obtain from \rf{A8} that
$$
L\gamma(x)\leq\varphi'(r(x))Lr(x)=\frac{1}{G(r(x)+1)}Lr(x)
\leq \mathrm{tr}P(x)\frac{e^{\int_0^{r(x)}G(s)ds}}{e^{\int_1^{r(x)}G(s)ds}-1}.
$$
Since $G\notin L^1(+\infty)$ we have
\be\label{A12}
\sup_{t\geq 2}\frac{e^{\int_0^{t}G(s)ds}}{e^{\int_1^{t}G(s)ds}-1}=\Lambda<+\infty.
\ee
We deduce that
$L\gamma(x)\leq \mathrm{tr}P(x)\Lambda$, i.e.,
\be\label{A14}
q(x)L\gamma(x)\leq \Lambda\;\;\mbox{on}\;\; D_o\cap(M^n\setminus \bar{B}_2).
\ee

Let $u\in\mathcal{C}^2(M)$ with $u^*=\sup_{M}u < +\infty$.
For a fixed $\eta>0$, consider
$$
A_\eta=\{ x\in M^n : u(x)>u^*-\eta \}
$$
and
$$
B_\eta=\{x\in A_\eta : \|\nabla u(x)\|<\eta \}.
$$
Since $M^n$ is complete, we have  from the  Ekeland quasi-minimum principle
(cf.\ \cite{Ek}) that $B_\eta\neq\emptyset$. All we have to show is that
\be\label{claim}
\inf_{B_\eta}\{q(x)Lu(x)\}\leq 0
\ee
since this is equivalent to the claim of the theorem.
To prove (\ref{claim})  we reason by contradiction. In fact, assume that
\be\label{A144}
q(x)Lu(x)\geq\sigma_0>0 \text{ on } B_\eta.
\ee
First observe that $u^*$ cannot be attained at a point $x_0\in M^n$, for otherwise
$x_0\in B_{\eta}$ but, since $P$ is positive semi-definite, then $q(x_0)Lu(x_0)\leq0$
thus contradicting (\ref{A144}). Set
$$
\Omega_t=\left\{x\in M^n:\gamma(x)>t\right\}.
$$
Then $\Omega^c_t=M^n\setminus \Omega_t$ is closed and hence compact by \rf{A15}. Define
$$
u^*_t=\max_{x\in\Omega^c_t}u(x).
$$

Since $u^*$ is not attained in $M^n$ and $\left\{\Omega^c_t\right\}$ is a
nested family exhausting $M^n$, there is a divergent
sequence $\left\{t_j\right\}_{j\in \enne}\subset[0,+\infty)$ such that
\be\label{A16}
u^*_{t_j}\rightarrow u^* \quad \hbox{as $j\rightarrow+\infty$},
\ee
and $T_1>0$ sufficiently large such that
$u^*_{T_1}>u^*-\eta/2$
and $\Omega_{T_1}\subset M^n\setminus \bar{B}_2$.
In particular, \rf{A14} holds on $D_o\cap\Omega_{T_1}$.
Choose $\alpha$ such that $u_{T_1}^*<\alpha<u^*$. Because of
(\ref{A16}) we can find $j$ sufficiently large such that
$T_2=t_j>T_1$ and $u^*_{T_2}>\alpha$. Then,
we select $\delta>0$ small enough so that
\be\label{A18}
\alpha+\delta<u^*_{T_2}.
\ee

For $\sigma>0$ define
$$
\gamma_{\sigma}(x)=\alpha+\sigma(\gamma(x)-T_1).
$$
Then, we have
$$
\gamma_{\sigma}(x)=\alpha \quad \hbox{for $x\in\partial\Omega_{T_1}$}
$$
and from \rf{A14} for $\sigma$ sufficiently small that
\be\label{A19}
q(x)L\gamma_{\sigma}(x)=\sigma q(x)L\gamma(x)\leq\sigma\Lambda<\sigma_0
\quad \hbox{on $D_o\cap\Omega_{T_1}$}.
\ee
On $\Omega_{T_1}\setminus\Omega_{T_2}$, we have
$$
\alpha\leq\gamma_{\sigma}(x)\leq\alpha+\sigma(T_2-T_1).
$$
Thus, choosing $\sigma>0$ sufficiently small so that
\be\label{A21}
\sigma(T_2-T_1)<\delta,
\ee
we obtain
$$
\alpha\leq\gamma_{\sigma}(x)<\alpha+\delta
\quad \hbox{on $\quad \Omega_{T_1}\setminus\Omega_{T_2}$}.
$$
For $x\in\partial\Omega_{T_1}$ we have that
$\gamma_{\sigma}(x)=\alpha>u^*_{T_1}\geq u(x)$.  Hence,
\be\label{A22}
(u-\gamma_{\sigma})(x)<0 \quad \hbox{on $\partial\Omega_{T_1}$}.
\ee
Let $\bar{x}\in\Omega_{T_1}\setminus\Omega_{T_2}$ be such that
$u(\bar{x})=u^*_{T_2}>\alpha+\delta$. Then (\ref{A18}) and \rf{A21}
yield
$$
(u-\gamma_{\sigma})(\bar{x})
\geq u^*_{T_2}-\alpha-\sigma(T_2-T_1)>u^*_{T_2}-\alpha-\delta>0.
$$
Moreover, we have from \rf{A15} and  $u^*<+\infty$ for $T_3>T_2$ sufficiently large that
\be\label{A23}
(u-\gamma_{\sigma})(x)<0 \quad \hbox{on $\Omega_{T_3}$}
\ee
Therefore,
$$
m=\sup_{x\in\bar{\Omega}_{T_1}}(u-\gamma_\sigma)(x)>0
$$
is, in fact,  a maximum attained at a point $z_0$ in the compact set
$\bar{\Omega}_{T_1}\setminus\Omega_{T_3}$.

 From \rf{A22} we know that  $\gamma(z_0)>T_1$. Thus,  we have
$$
u(z_0)=\gamma_{\sigma}(z_0)+m>\gamma_{\sigma}(z_0)>\alpha>u^*_{T_1}>u^*-\frac{\eta}{2},
$$
and hence $z_0\in A_{\eta}\cap\Omega_{T_1}$.  Next, we have to distinguish two
cases, according to $z_0\in D_o$ or $z_0\notin D_o$.

If $z_0\in D_o$, since $z_0$ is a maximum for $u-\gamma_\sigma$, we have
$\nabla(u-\gamma_\sigma)(z_0)=0$. Using this fact, we have that
$z_0\in B_\eta$ since
$$
\|\nabla u(z_0)\|=\|\nabla\gamma_\sigma(z_0)\|
=\sigma\varphi'(r(z_0))\|\nabla r(z_0)\|=\frac{\sigma}{G(r(z_0)+1)}
\leq\frac{\sigma}{G(1)}<\eta
$$
up to choosing $\sigma$ sufficiently small. Since $P$ is positive semi-definite
and $z_0$ is a maximum for $u-\gamma_\sigma$, we have
$Lu(z_0)\leq L\gamma_\sigma(z_0)$, and this jointly with \rf{A14} yields
$$
0<\sigma_0\leq q(z_0)Lu(z_0)\leq q(z_0)L\gamma_{\sigma}(z_0)<\sigma_0,
$$
which is a contradiction and concludes the proof for this case.

In the case $z_0\notin D_o$ we reason as follows. Fix $0<\varepsilon<1$
sufficiently small so that for the minimizing geodesic
$\varsigma$ parametrized by arclength and joining $o$ with $z_0$,
the point $o_\varepsilon=\varsigma(\varepsilon)\neq z_0$
and $z_0\notin\mathrm{cut}(o_\varepsilon)$.
Hence, the function $r_\varepsilon(x)=\mathrm{dist}(o_\varepsilon,x)$
is $\mathcal{C}^2$ in a neighborhood of $z_0$. By the triangle inequality
\be\label{A24.1}
r(x)\leq r_\varepsilon(x)+\varepsilon,
\ee
equality holding at $z_0$. With $\varphi$ defined in \rf{A10}  set
$$
\gamma^\varepsilon(x)=\varphi(r_\varepsilon(x)+\varepsilon).
$$
Since $\varphi$ is increasing
\be\label{A25}
\gamma(x)=\varphi(r(x))\leq\varphi(r_\varepsilon(x)+\varepsilon)=\gamma^\varepsilon(x).
\ee
and
\be\label{A26}
\gamma(z_0)=\gamma^\varepsilon(z_0).
\ee
Next consider the function
$$
\gamma_\sigma^\varepsilon(x)=\alpha+\sigma(\gamma^\varepsilon(x)-T_1)\; (\geq \gamma_\sigma(x)).
$$
Because of \rf{A25} and \rf{A26} we have in a neighborhood of $z_0$ that
$$
u(x)-\gamma_\sigma^\varepsilon(x)\leq u(x)-\gamma_\sigma(x)\leq m
$$
and
$$
u(z_0)-\gamma_\sigma^\varepsilon(z_0)=u(z_0)-\gamma_\sigma(z_0)=m.
$$
Hence $z_0$ is also a local maximum for $u(x)-\gamma_\sigma^\varepsilon(x)$. Thus,
\be\label{A27}
\nabla u(z_0)=\nabla\gamma_\sigma^\varepsilon(z_0)
\ee
and
\be\label{A28}
Lu(z_0)\leq L\gamma_\sigma^\varepsilon(z_0).
\ee
 From \rf{A27} we deduce
\bea
\|\nabla u(z_0)\|\!\!\!& = &\!\!\! \sigma\|\nabla\gamma^\varepsilon(z_0)\|=\sigma\varphi'(r_\varepsilon(z_0)
+\varepsilon)\|\nabla r_\varepsilon(z_0)\|\\
\!\!\!& = &\!\!\!\frac{\sigma}{G(r_\varepsilon(z_0)+\varepsilon+1)}\leq\frac{\sigma}{G(1)}<\eta.
\eea
Since we already knew that $z_0\in A_\eta$, we conclude that $z_0\in B_\eta$.
Now we analyze \rf{A28}. Because of \rf{A1}, \rf{A24.1} and $G'\geq 0$ we
have
$$
K_M(x)\geq -G^2(r(x))\geq -G^2(r_\varepsilon(x)+\varepsilon).
$$
Set $G_\varepsilon(t)=G(t+\varepsilon)$
and consider the Cauchy problem \rf{A4} with $G_\varepsilon$ instead of $G$.
Again by the Hessian comparison theorem,  on $D_{o_\varepsilon}$ we have
$$
Lr_\varepsilon(x)\leq
\mathrm{tr}P(x)\frac{\psi'_\varepsilon(r_\varepsilon(x))}{\psi_\varepsilon(r_\varepsilon(x))}
\;\;\;\mbox{where}\;\;\;
\psi_\varepsilon(t)=\frac{1}{G_\varepsilon(0)}\left(e^{\int_0^tG_\varepsilon(s)ds}-1\right).
$$
Observing that $z_0\in D_{o_\varepsilon}$, we obtain using \rf{A12} that
\bea
L\gamma^\varepsilon(z_0) \!\!\!& \leq &\!\!\!\varphi'(r_\varepsilon(z_0)+\varepsilon)Lr_\varepsilon(z_0)
=\frac{1}{G(r_\varepsilon(z_0)+\varepsilon+1)}Lr_\varepsilon(z_0)\\
\!\!\!& = &\!\!\! \frac{1}{G(r(z_0)+1)}Lr_\varepsilon(z_0)
\leq\frac{\mathrm{tr}P(z_0)}{G(r(z_0)+1)}\frac{\psi'_\varepsilon(r_\varepsilon(z_0))}{\psi_\varepsilon(r_\varepsilon(z_0))}\\
\!\!\!& = &\!\!\!\frac{\mathrm{tr}P(z_0)}{G(r(z_0)+1)}G_\varepsilon(r_\varepsilon(z_0))\frac{e^{\int_0^{r_\varepsilon(z_0)}G(s+\varepsilon)ds}}{e^{\int_0^{r_\varepsilon(z_0)}G(s+\varepsilon)ds}-1}\\
\!\!\!& = &\!\!\!\mathrm{tr}P(z_0)\frac{G(r_\varepsilon(z_0)+\varepsilon)}{G(r(z_0)+1)}\frac{e^{\int_\varepsilon^{r_\varepsilon(z_0)+\varepsilon}G(s)ds}}{e^{\int_\varepsilon^{r_\varepsilon(z_0)+\varepsilon}G(s)ds}-1}\\
\!\!\!& = &\!\!\!\mathrm{tr}P(z_0)\frac{G(r(z_0))}{G(r(z_0)+1)}
\frac{e^{\int_\varepsilon^{r(z_0)}G(s)ds}}{e^{\int_\varepsilon^{r(z_0)}G(s)ds}-1}\\
\!\!\!&\leq&\!\!\!\mathrm{tr}P(z_0)\frac{e^{\int_0^{r(z_0)}G(s)ds}}{e^{\int_1^{r(z_0)}G(s)ds}-1}
\leq\mathrm{tr}P(z_0)\Lambda.
\eea
Thus,
$$
L\gamma_\sigma^\varepsilon(z_0)=\sigma L\gamma^\varepsilon(z_0)
\leq\mathrm{tr}P(z_0)\sigma\Lambda<\mathrm{tr}P(z_0)\sigma_0.
$$
 From \rf{A14} and \rf{A28} we  deduce that
$$
0<\sigma_0\leq q(z_0)Lu(z_0)\leq q(z_0)L\gamma^\varepsilon_{\sigma}(z_0)
\leq\sigma\Lambda<\sigma_0,
$$
and this is a contradiction.\qed

\section{The Riemannian case}

Let $f\colon M^n\to\bar{M}^{n+1}$ denote an isometric immersion
between Riemannian manifolds. Assume that the hypersurface $f$ is
two-sided, that is, there exists a globally defined unit normal vector field
$N$. Denote by $A=A_N$ the second fundamental form of $f$ for the given
orientation Then, the \emph{$k$-mean curvature} $H_k$ is given by
$$
{n \choose k}H_k= S_k,\;\;0\leq k\leq n,
$$
where $S_0=1$ and $S_k$ for $k\geq 1$ is the $k$-symmetric elementary function on
the principal curvatures of $f$. In particular, when $k=1$ then $H_1=H$ is the mean
curvature of $f$. Moreover, for $k$ even the sign of $S_k$
(and hence $H_k$) does not depend on the chosen orientation.

The Newton tensors $P_k\colon TM\rightarrow TM$, $0\leq k\leq n$, arising from $A$ are
defined inductively by $P_0=I$ and $P_k=S_k I-AP_{k-1}$.
Then,
\be\label{traces}
\mathrm{Tr}P_k=(n-k)S_k=c_kH_k\;\;\;\; \mbox{and}\;\;\;\;
\mathrm{Tr}AP_k=(k+1)S_{k+1}=c_{k}H_{k+1}
\ee
where
$c_k=(n-k){n \choose k}=(k+1){n \choose k\!+\!1}$.

The second order differential operators
$L_k\colon\mathcal{C}^{\infty}(M)\rightarrow\mathcal{C}^{\infty}(M)$
arise from normal variations of $P_{k+1}$ and are given  by
$$
L_k=\mathrm{Tr}\,(P_k \circ \hess).
$$
Then, the operator $L_k$ is semi-elliptic (respectively, elliptic) if and only
if $P_k$ is positive semi-definite (respectively, positive definite).

Let $B_{\bar{M}}(r)$ denote the geodesic ball  with radius $r$
centered  at a  reference point $o\in\bar{M}^{n+1}$. In the
sequel, we assume that the radial sectional curvatures in
$B_{\bar{M}}(r)$ along the geodesics issuing from $o$  are bounded
as $K_{\bar{M}}^{\mathrm{rad}}\leq b$
for some constant $b\in\R$, and that $0<r<\min\{{\rm
inj}_{\bar{M}}(o), \pi/2\sqrt{b}\}$ where ${\rm inj}_{\bar{M}}(o)$
is the injectivity radius at $o$ and $\pi/2\sqrt{b}$ is replaced by
$+\infty$ if $b\leq 0$.

It is a standard fact that if  $\bar{M}^{n+1}$ has constant sectional
curvature $b$, then the mean curvature of the geodesic sphere
$S_{\bar{M}}(r)=\partial B_{\bar{M}}(r)$ is
\be\label{cbr}
C_b(r)
=\left\{\begin{array}{lll}
\sqrt{b}\cot(\sqrt{b}\, r) & \mathrm{if} & b >0,\\
1/r & \mathrm{if} & b =0,\\
\sqrt{-b}\coth(\sqrt{-b}\, r) & \mathrm{if}  & b <0.
\end{array}\right.
\ee
The following classical Hessian comparison result plays an important role in
the proof of our results.
\begin{lemma}\po\label{lemmariemannian}
Let $\bar{M}$ be a Riemannian manifold with a fixed reference point
$o\in\bar{M}$ and let $\rho(x)$ be the distance function to $x$.
Let $x\in\bar{M}$ be inside a geodesic ball $B_{\bar{M}}(r)$ as above with
$K_{\bar{M}}^{\mathrm{rad}}\leq b$.
Then, for any vector $X\in T_xM$ we have
$$
\Hess{\rho}(X,X)\geq C_b(\rho(x))(\|X\|^2-\<X,\bar{\nabla}\rho(x)\>^2)
$$
where $\Hess{\rho}$ stands for the Hessian of $\rho$.
\end{lemma}

In the following result, it is convenient to think that
$S_{\bar{M}}(r)$ is the smallest possible geodesic sphere centered at  $o$
enclosing the hypersurface.

\begin{theorem}\po\label{main} Let $f\colon M^n\to\bar{M}^{n+1}$ be a
two-sided isometric immersion between complete manifolds where
$M^n$ satisfies condition (\ref{A1}). Assume that $P_k$ is
positive semi-definite for some $0\leq k\leq n-1$ and that $\mathrm{tr}P_k>0$ on $M^n$.
If $f(M)\subset B_{\bar{M}}(r)$
for a geodesic ball $B_{\bar{M}}(r)$ as above, then
\be\label{ineq}
\sup_M\left(\frac{|H_{k+1}|}{H_k}\right)\geq C_b(r).
\ee
Moreover, if  $P_k$ is positive definite and there exists a point  $p_0\in
M^n$ such that $f(p_0)\in S_{\bar{M}}(r)$ then equality in (\ref{ineq}) implies
$M^n=S_{\bar{M}}(r)$.
\end{theorem}

In particular, we have the following consequence.

\begin{corollary}\po\label{coromain} Let $f\colon M^n\to\bar{M}^{n+1}$ be as above.
Assume that $P_k$ is positive semi-definite for some $0\leq k\leq
n-1$. If $f(M)\subset B_{\bar{M}}(r)$ for a geodesic ball
$B_{\bar{M}}(r)$ as above, then
\be\label{ineqcoro}
\sup_M|H_{k+1}|\geq C_b(r)\,\inf_M H_k.
\ee
\end{corollary}

For the proof of Corollary \ref{coromain} we first observe that
(\ref{ineqcoro}) holds trivially if \mbox{$\inf_M H_k=0$}. For $\inf_M H_k>0$, we have
that $P_k\neq0$ everywhere and the result follows directly from Theorem
\ref{main} since (\ref{ineqcoro}) is
weaker than (\ref{ineq}).
\vspace{2ex}

\noindent{\it Proof of Theorem \ref{main}:} We denote
by $\rho\colon \bar{M}^{n+1}\to\R$ the distance function to the
reference point $o$ and set $u=\rho\circ f$. Along $M^n$ we have
$$
\bar\nabla\rho=\nabla u + \<\bar\nabla\rho,N\>N
$$
where $N$ is a unit global normal vector field to $f$.  An easy
computation gives
$$
\Hess{u}(X,Y)=\Hess{\rho}(X,Y) + \<\bar\nabla\rho,N\>\<AX,Y\>
$$
where we denoted  $A=A_N$.
\newpage

Let $e_1,\ldots,e_n$ be an orthonormal basis of principal directions at a
point of $M^n$. We obtain using (\ref{traces}) that
\bea
L_ku\!\!\!&=&\!\!\!\sum_{i=1}^n\Hess{u}(e_i,P_ke_i)
=\sum_{i=1}^n\Hess{\rho}(e_i,P_ke_i) + \<\nap\rho,N\>\mathrm{Tr} AP_k\\
\!\!\!&=&\!\!\!\sum_{i=1}^n\Hess{\rho}(e_i,P_ke_i)+ c_kH_{k+1}\<\nap\rho,N\>.
\eea
By assumption, we have
$$
P_ke_i=\mu_ie_i\;\;\;\mbox{with}\;\;\;\mu_i\geq 0.
$$
Using the Hessian comparison theorem, we obtain
\bea
\Hess{\rho}(e_i,P_ke_i)\!\!\!&=&\!\!\!\mu_i\Hess{\rho}(e_i,e_i)\\
\!\!\!&\geq&\!\!\!\mu_iC_b(u)(1-\<\nabla u,e_i\>^2)\\
\!\!\!&=&\!\!\!C_b(u)(\mu_i-\<\nabla u,e_i\>\<P_k\nabla u,e_i\>).
\eea
Using (\ref{traces}) we have
\bea
\sum_{i=1}^n\Hess{\rho}(e_i,P_ke_i)\!\!\!&\geq&\!\!\! C_b(u)(\mathrm{Tr}P_k
-\<\nabla u,P_k\nabla u)\>\\
\!\!\!&=&\!\!\!C_b(u)(c_kH_k-\<\nabla u,P_k\nabla u\>).
\eea
Therefore,
\be\label{one}
L_ku\geq C_b(u)(c_kH_k-\<\nabla u,P_k\nabla u\>)+c_kH_{k+1}\<\nap\rho,N\>.
\ee

Consider the function
$$
\phi_b(t)
=\left\{\begin{array}{lll}
1-\cos(\sqrt{b}\,t) & \mathrm{if} & b >0,\\
t^2 & \mathrm{if} & b =0,\\
\coth(\sqrt{-b}\,t) & \mathrm{if}  & b <0.
\end{array}\right.
$$
Then $\phi_b'(t)>0$ if $t>0$ and
\be\label{two}
\phi_b''(t) - C_b(t)\phi_b'(t)=0.
\ee
We have using (\ref{two}) that
\bea
L_k\phi_b(u)\!\!\!&=&\!\!\!\phi_b''(u)\<\nabla u, P_k\nabla u\>+\phi_b'(u)L_ku\\
\!\!\!&=&\!\!\!\phi_b'(u)(C_b(u)\<\nabla u, P_k\nabla u\>+L_ku).
\eea
It follows from (\ref{one}) that
$$
L_k\phi_b(u)\geq c_k\phi_b'(u)(C_b(u)H_k+\<\nap\rho,N\>H_{k+1}).
$$
Hence,
$$
L_k\phi_b(u)\geq c_k\phi_b'(u)\left(C_b(u)H_k-|H_{k+1}|\right).
$$
Since $\sup_M\phi_b(u)\leq \phi_b(r)<+\infty$, it follows from Theorem \ref{MAIN}
that there exists a sequence of points $\left\{p_{j}\right\}_{j\in \enne}$ in $M^n$
such that
$$
\phi_b(u(p_j))>\sup_M\phi_b(u)-\frac{1}{j}\;\;\;
\mbox{and} \;\;\;
\frac{1}{c_kH_k(p_j)}L_k\phi_b(u)(p_j)< \frac{1}{j}.
$$
It follows from the first inequality that
\be\label{u}
\lim_{j\rightarrow\infty}u(p_j)=u^*=\sup_Mu
\ee
since $\sup_M\phi_b(u)=\phi_b(\sup_Mu)$. Therefore,
\bea
\frac{1}{j}>\frac{1}{c_kH_k(p_j)}L_k\phi_b(u)(p_j) \!\!\!& \geq &\!\!\!
\phi_b'(u(p_j))\left(C_b(u(p_j))-\frac{|H_{k+1}|}{H_k}(p_j)\right)\\
\!\!\!& \geq &\!\!\! \phi_b'(u(p_j))\left(C_b(r)-\sup_M\left(\frac{|H_{k+1}|}{H_k}\right)\right)
\eea
since $C_b(u(p_j)\geq C_b(r)$. Taking $j\to +\infty$ and using (\ref{u}) we conclude that
$$
C_b(r)-\sup_M\left(\frac{|H_{k+1}|}{H_k}\right)\leq 0.
$$

For the proof of the second statement, first observe that equality in (\ref{ineq})
yields $L_k\phi_b(u)\geq 0$. Since $\phi_b(u)\leq\phi_b(r)<+\infty$, it follows
from the maximum principle for the elliptic operator $L_k$ that $\phi_b(u)$ is
constant, and hence $u$ is constant.\qed

\begin{remark}\po {\em Notice that the conclusion $(ii)$ in Theorem \ref{maxprinc}
has not been used in the proof of Theorem \ref{main}. In this situation, the usual
terminology is that we only need a \emph{weak} Omori-Yau maximum principle
for trace operators. It turns out that for spacelike hypersurfaces
in Lorentzian ambient spaces this is not longer the case.}
\end{remark}

In the sequel, we replace some  assumptions in Theorem \ref{main} by
simpler ones and of a geometric nature. This, of course, is the case of
Theorem~\ref{main3}
in the Introduction.  But first we considered the special case of $H_2$.
The short proofs given next are mostly taken from \cite{AIR} and are included
for the sake of completeness.

\begin{corollary}\po\label{main2}\!\!Let $f\colon M^n\to\bar{M}^{n+1}$ be an isometric
immersion into a \mbox{complete} Riemannian manifold. Assume that $M^n$ is complete with
sectional curvature $K_M\geq K>-\infty$. If $H_2>0$ and
$f(M)\subset B_{\bar{M}}(r)$ for~a geodesic ball $B_{\bar{M}}(r)$ as above,
then
\be\label{ineq2}
\sup_M\sqrt{H_2}\geq\sup_M\left(\frac{H_2}{H}\right)\geq C_b(r).
\ee
If  there exists a point $p_0\in M^n$ such that $f(p_0)\in S_{\bar M}(r)$
and it holds that $\sup_M(H_2/H)=C_b(r)$,  then $M^n=S_{\bar{M}}(r)$.
\end{corollary}

\proof   In term of the principal curvatures $\lambda_1,\ldots,\lambda_n$ of $f$
we have that
$$
n^2H^2=\sum_{j=1}^n\lambda_j^2+n(n-1)H_2>\lambda_i^2.
$$
In particular, the immersion is two-sided since $H^2>0$. Moreover, we have that the
eigenvalues of $P_1$ satisfy $\mu_j=nH-\lambda_j>0$ for any $j$ (see Lemma~3.10 in \cite{E})
and therefore $L_1$ is elliptic. Then, the second inequality and the characterization of equality
follows from Theorem  \ref{main}. For the first inequality, just observe that $H^2-H_2\geq 0$
yields $H_2/H\leq \sqrt{H_2}$.
\qed

\begin{remark}\po{\em
If the ambient space has constant sectional curvature $b$, then the  normalized
scalar curvature $s$ of $M^n$ is related to $H_2$ by $s=b+H_2$. In this case
inequality (\ref{ineq2}) gives
$$
\sup_Ms\geq b+C_b(r)\,\inf_M H.
$$
}\end{remark}

\noindent{\it Proof of Theorem \ref{main3}:} The existence of an elliptic point implies
that $H_{k+1}$ is positive at that
point, and hence on $M^n$. The well-known Garding inequalities yield,  for the
appropriate orientation, that
\be\label{garding}
H_1\geq H_2^{1/2}\geq\cdots\geq H_{k}^{1/k}\geq H_{k+1}^{1/(k+1)}>0.
\ee
Thus, the immersion is two-sided and $H_1>0$. Moreover,
since $M^n$ has an elliptic point and $H_{k+1}\neq 0$ on $M^n$,
from the proof of \cite[Proposition 3.2]{BC} we have that the operators
$L_{j}$ are elliptic for any $1\leq j\leq k$.
Then, the second inequality and the characterization of the
equality case follows from Theorem~\ref{main}. For the first inequality observe
that $H_{j+1}/H_j\leq\sqrt[j+1]{H_{j+1}}$
follows from (\ref{garding}).\qed

\section{The Lorentzian case}

Let $f\colon M^n\to\bar{M}^{n+1}$ be a spacelike hypersurface
isometrically immersed into a spacetime. Since $\bar{M}^{n+1}$ is
time-oriented, there exists a unique globally defined
future-directed timelike normal unit vector $N$.
We refer to $N$ as the future-directed Gauss map of $M^n$ and
denote by $A=A_N$ the second fundamental form of the
hypersurface.

For spacelike hypersurfaces, the $k$-mean curvature $H_k$ is defined
by
$$
{n \choose k}H_k= (-1)^kS_k,\;\;0\leq k\leq n,
$$
where $S_0=1$ and $S_k$ for $k\geq 1$ is the $k$-symmetric
elementary function on the principal curvatures of $f$. The choice
of the sign $(-1)^k$ in the definition is to have the
mean curvature vector given by $\overrightarrow{H}=HN$.
Therefore, $H(p)>0$ at $p\in M^n$ if and only if
$\overrightarrow{H}(p)$ is future-directed. Clearly, when $k$ is
even the sign of $H_k$ does not depend on the chosen Gauss map.

For spacelike hypersurfaces, the Newton tensors $P_k\colon
TM\rightarrow TM$ are defined inductively by $P_0=I$ and
$P_k=(-1)^kS_k I+AP_{k-1}$, $1\leq k\leq n$.~Then,
\be\label{tracesL}
\mathrm{Tr}P_k=c_kH_k\;\;\;\; \mbox{and}\;\;\;\;\mathrm{Tr}AP_k=-c_{k}H_{k+1}.
\ee

Let $o\in\bar{M}^{n+1}$ be a reference point and $\rho\colon
\bar{M}^{n+1}\rightarrow [0,+\infty]$ the Lorentzian distance
from $o$. It is well known that the Lorentzian distance function may
fail to be continuous and even finite valued. Thus,  to
guarantee smoothness we need to restrict $\rho$
to certain special subsets of $\bar{M}^{n+1}$. Following
\cite{EGK} (see also \cite{AHP}) we denote by
$\mathcal{I}^{+}(o)\subset\bar{M}^{n+1}$ the diffeomorphic image of
$\textrm{int}(\tilde{\mathcal{I}}^{+}(o))$ under the exponential map
at $o$.  Here,
$$
\tilde{\mathcal{I}}^{+}(o)=\{ tv\in T_o\bar{M}: \mbox{$v$ future-directed unit vector}
\mbox{ and } 0<t<s_o(v)\}
$$
where
$$
s_o(v)=\sup\{ t\geq 0 : \rho(\gamma_v(t))=t=L(\gamma_v|_{[0,t]}) \}.
$$
It turns out that $\mathcal{I}^{+}(o)$ is the largest natural open
subset of $\bar{M}^{n+1}$ on which $\rho$ is smooth and that
$\bar\nabla \rho$ is a {past-directed timelike (geodesic) unit}
vector field on $\mathcal{I}^{+}(o)$. We refer  to \cite{AHP},
\cite{EGK} and references therein for further details about the
Lorentzian distance function.

For $b\in\R{}$, we consider the function
$\widehat{C}_b(t)=C_{-b}(t)$.
We point out that when $\mathcal{I}^{+}(o)\neq\emptyset$,
then $\widehat{C}_b(r)$ is the future mean curvature of the level set
$$
\Sigma_b(r)=\{ x \in\mathcal{I}^{+}(o) : \rho(x)=r \}
$$
in a Lorentzian space form $\bar{M}^{n+1}_b$ with constant sectional
curvature $b$. The following Hessian comparison result  plays an
important role in the proof of our  results.
\begin{lemma}\po\label{lemmalorentz}
Assume that $\mathcal{I}^{+}(o)\neq\emptyset$ for a reference point
$o\in\bar{M}^{n+1}$. Let $x\in\mathcal{I}^{+}(o)$ and assume that the
radial sectional curvatures of $\bar{M}^{n+1}$ along the radial
future geodesic from $o$ to $x$ are bounded as
$K^{\mathrm{rad}}_{\bar{M}}\leq b$ (respectively, $K^{\mathrm{rad}}_{\bar{M}}\geq b$)
for some constant $b$, with $\rho(x)<\pi/\sqrt {-b}$ if $b<0$. Then, for
any spacelike vector $X\in T_xM$ we have
$$
\Hess{\rho}(X,X)\geq
-\widehat{C}_b(\rho(x))(\|X\|^2+\<X,\bar{\nabla}\rho(x)\>^2) \quad
\text{(respectively, }  \leq  \text{)}
$$
where $\Hess{\rho}$ stands for the Lorentzian Hessian of $\rho$.
\end{lemma}

The proof of Lemma \ref{lemmalorentz} follows easily from the proofs
of Lemma 3.1 and Lemma 3.2 in \cite{AHP} by observing that the
assumption $K_{\bar{M}}\leq b$ (respectively, $K_{\bar{M}}\geq b$)
for all timelike planes in $\bar{M}^{n+1}$ in those results is now
needed only for the radial sectional curvatures along the radial
future geodesic starting at $o$. Observe also that
$$
\Hess{\rho}(X,X)=\Hess{\rho}(X^*,X^*)
$$
where $X=X^*-\<X,\bar{\nabla}\rho(x)\>\bar\nabla\rho(x)$ with
$\<X^*,\bar{\nabla}\rho(x)\>=0$ and
$$
\|X^*\|^2= \|X\|^2+\<X,\bar{\nabla}\rho(x)\>^2.
$$
For details, see \cite[Section 3]{AHP}.

For a given reference point $o\in\bar{M}^{n+1}$ and $r>0$,
let $B^{+}(o,r)$ denote the future inner ball of radius $r$, namely,
$$
B^{+}(o,r)=\{ x\in I^{+}(o) : \rho(x)<r\},
$$
where $I^{+}(o)$ is the chronological future of $o$, i.e., the
set of points \mbox{$x\in\bar{M}^{n+1}$} for which there exists a
future-directed timelike curve from $o$ to $x$. Now we are ready to
state our first result in this section.

\begin{theorem}\po\label{lorentz}
Let $f\colon M^n\to\bar{M}^{n+1}$ be a spacelike
hypersurface immersed into a spacetime, where $M^n$ is complete  and
satisfies condition (\ref{A1}). Assume that $P_k$ is positive
semi-definite for some $0\leq k\leq n-1$  and that
$\mathrm{tr}P_k>0$ on $M^n$. If
$f(M)\subset\mathcal{I}^{+}(o)\cap B^{+}(o,r)$ for a reference point
$o\in\bar{M}^{n+1}$ and the radial sectional curvatures along the
radial future geodesics issuing from $o$ are bounded as
$K^{\mathrm{rad}}_{\bar{M}}\leq b$ on $\mathcal{I}^{+}(o)\cap
B^{+}(o,r)$ for some constant $b$ (with $r<\pi/2\sqrt{-b}$ if
$b<0$), then
\be\label{ineqlorentz1}
\inf_M\left(\frac{H_{k+1}}{H_k}\right)\leq \widehat{C}_b(u^{*}), \ee
where $u^{*}=\sup_Mu$.
\end{theorem}

\proof
We  set $u=\rho\circ f$. Along $M^n$, we have
$$
\bar\nabla\rho=\nabla u - \<\bar\nabla\rho,N\>N= \nabla u-\sqrt{1+\|\nabla u\|^2}N.
$$
Then,
$$
\Hess{u}(X,Y)=\Hess{\rho}(X,Y)-\sqrt{1+\|\nabla u\|^2}\<AX,Y\>.
$$
A similar computation as in the Riemannian case yields
$$
L_ku=\sum_{i=1}^n\Hess{\rho}(e_i,P_ke_i)+ c_kH_{k+1}\sqrt{1+\|\nabla u\|^2},
$$
where $e_1,\ldots,e_n$ is an orthonormal basis of principal directions at a point of $M^n$.
Using Lemma \ref{lemmalorentz} and reasoning as we did  to conclude (\ref{one}), we obtain
$$
L_ku\geq-\widehat{C}_b(u)(c_kH_k+\<\nabla u,P_k\nabla u
\>)+c_kH_{k+1}\sqrt{1+\|\nabla u\|^2},
$$
where the restriction $u<\pi/2\sqrt{-b}$ if $b<0$ is necessary to have $\widehat{C}_b(u)>0$.

Since $\sup_Mu=u^*<+\infty$, by Theorem \ref{MAIN} there exists
a sequence  $\left\{p_{j}\right\}_{j\in \enne}$ in  $M^n$ such that
$$
u(p_j)>u^*-\frac{1}{j},\;\;
\|\nabla u(p_j)\|<\frac{1}{j}\;\; \mbox{and} \;\;
\frac{1}{c_kH_k(p_j)}L_ku(p_j)< \frac{1}{j}.
$$
In particular, $\lim_{j\rightarrow\infty}u(p_j)
=u^*$ and $\lim_{j\rightarrow\infty}\|\nabla u(p_j)\|=0$. Thus,
\begin{eqnarray*}
\frac{1}{j} \!\!\!& > &\!\!\! \frac{1}{c_kH_k(p_j)}L_ku(p_j)\\
\!\!\!& \geq &\!\!\! -\widehat{C}_b(u(p_j))\left(1
+\frac{\<\nabla u,P_k\nabla u \>}{c_kH_k}(p_j)\right)
+\frac{H_{k+1}}{H_k}(p_j)\sqrt{1+\|\nabla u(p_j)\|^2}\\
\!\!\!& \geq &\!\!\! -\widehat{C}_b(u(p_j))\left(1+\|\nabla
u(p_j)\|^2\right)+\inf_M\left(\frac{H_{k+1}}{H_k}\right)\sqrt{1+\|\nabla
u(p_j)\|^2},
\end{eqnarray*}
since, being $P_k$ positive semi-definite, we have
$$
0\leq\<X,P_kX\>\leq\mathrm{Tr}P_k\|X\|^2=c_kH_k\|X\|^2
$$
for any $X\in TM$. Finally, taking $j\to +\infty$  we conclude that
(\ref{ineqlorentz1}) holds. \qed\vspace{1,5ex}

The following is the second main result in this section.

\begin{theorem}\po\label{lorentzB}
Let $f\colon M^n\to\bar{M}^{n+1}$ be a spacelike
hypersurface immersed into a spacetime, where $M^n$ is complete  and
satisfies condition (\ref{A1}). Assume that $P_k$ is positive
semi-definite for some $0\leq k\leq n-1$  and that
$\mathrm{tr}P_k>0$ on $M^n$. If $f(M)\subset\mathcal{I}^{+}(o)$
for a reference point $o\in\bar{M}^{n+1}$ and the radial sectional
curvatures along the radial future geodesics issuing from $o$ are
bounded as $K^{\mathrm{rad}}_{\bar{M}}\geq b$ on
$\mathcal{I}^{+}(o)$ for some constant $b$ (with $r<\pi/2\sqrt{-b}$
if $b<0$), then
\be\label{ineqlorentz1B}
\sup_M\left(\frac{H_{k+1}}{H_k}\right)\geq \widehat{C}_b(u_{*}), \ee
where $u_{*}=\inf_Mu$. In particular,  if $u_*=0$ then
$\sup_M\left(H_{k+1}/H_k\right)=+\infty$.
\end{theorem}

\proof
We proceed as in the proof of Theorem \ref{lorentz}, by observing that
in this case Lemma \ref{lemmalorentz} yields
\begin{eqnarray*}
L_ku \!\!\!& \leq &\!\!\! -\widehat{C}_b(u)(c_kH_k+\<\nabla u,P_k\nabla u \>)
+c_kH_{k+1}\sqrt{1+\|\nabla u\|^2}\\
\!\!\! & \leq &\!\!\!
-c_kH_k\widehat{C}_b(u)+c_kH_{k+1}\sqrt{1+\|\nabla u\|^2}.
\end{eqnarray*}
Since $\inf_Mu=u_*\geq 0$, by Theorem \ref{MAIN} there is a
sequence  $\left\{p_{j}\right\}_{j\in \enne}\subset M^n$ such that
$$
u(p_j)<u_*+\frac{1}{j},\;\;
\|\nabla u(p_j)\|<\frac{1}{j}\;\; \mbox{and} \;\;\frac{1}{c_kH_k(p_j)}L_ku(p_j)>
-\frac{1}{j}.
$$
In particular, $\lim_{j\rightarrow\infty}u(p_j)=u_*$ and
$\lim_{j\rightarrow\infty}\|\nabla u(p_j)\|=0$. Thus
\begin{eqnarray*}
-\frac{1}{j} \!\!\!& < &\!\!\! \frac{1}{c_kH_k(p_j)}L_ku(p_j)\\
\!\!\!& \leq &\!\!\! -\widehat{C}_b(u(p_j))
+\frac{H_{k+1}}{H_k}(p_j)\sqrt{1+\|\nabla u(p_j)\|^2}\\
\!\!\!& \leq &\!\!\!
-\widehat{C}_b(u(p_j))+\sup_M\left(\frac{H_{k+1}}{H_k}\right)\sqrt{1+\|\nabla
u(p_j)\|^2},
\end{eqnarray*}
and we conclude taking $j\to +\infty$ that (\ref{ineqlorentz1B})
holds. The last assertion follows from (\ref{ineqlorentz1B}) and the
fact that $\lim_{t\rightarrow 0^+}\widehat{C}_b(t)=+\infty$.
\qed\vspace{1,5ex}

As a direct application of Theorem \ref{lorentzB} we get the
following result.

\begin{corollary}\po
Under the assumptions of Theorem \ref{lorentzB}, assume also that $H_{k+1}/H_k$
is
bounded from above on $M^n$.  Then,  there exists $\delta>0$ such that
$f(M)\subset O^+(o,\delta)$, where $O^+(o,\delta)$ denotes the
future outer ball of radius $\delta$ given by
$$
O^+(o,\delta)=\{ x\in I^+(o) : \rho(x)>\delta\}.
$$
\end{corollary}

\proof Simply observe that $\sup_M(H_{k+1}/H_k)<+\infty$ implies  $u_*>0$.\qed
\vspace{1,5ex}

If the ambient spacetime is a Lorentzian space
form, from  Theorem \ref{lorentz} and Theorem \ref{lorentzB} we obtain
the following consequence that extends Theorem~4.5 in \cite{AHP} to the case
of higher order mean curvatures.

\begin{corollary}\po\label{lorentzC}
Let $f\colon M^n\to\bar{M}^{n+1}$ be a
spacelike hypersurface immersed into a Lorentzian spacetime of constant
sectional curvature $b$, where $M^n$ is complete  and satisfies
condition (\ref{A1}). Assume that $P_k$ is positive semi-definite
for some $0\leq k\leq n-1$  and that $\mathrm{tr}P_k>0$ on $M^n$.
If $f(M)\subset\mathcal{I}^{+}(o)\cap B^+(o,r)$ for a reference
point $o\in\bar{M}$ (with $r<\pi/2\sqrt{-b}$ if $b<0$), then
\be\label{ineqlorentz1C}
\inf_M\left(\frac{H_{k+1}}{H_k}\right)\leq
\widehat{C}_b(u^{*}) \leq\widehat{C}_b(u_{*})\leq
\sup_M\left(\frac{H_{k+1}}{H_k}\right).
\ee
\end{corollary}

\vspace*{-1ex}

{\renewcommand{\baselinestretch}{1}
\hspace*{-20ex}\begin{tabbing}
\indent \=Luis J. Alias \hspace{25,5ex} Marcos Dajczer  \\
\>Departamento de Matematicas  \hspace{7,5ex} IMPA\\
\>Universidad de Murcia \hspace{15,4ex}
Estrada Dona Castorina, 110\\
\>  Campus de Espinardo
\hspace{17ex}22460-320 --- Rio de Janeiro ---RJ\\
\> E-30100 Espinardo, Murcia\hspace{12,6ex}Brazil \\
\> Spain\hspace{33,6ex} marcos@impa.br\\
\> ljalias@um.es
\end{tabbing}}
\vspace*{-1ex}
{\renewcommand{\baselinestretch}{1}
\hspace*{-20ex}\begin{tabbing}
\indent \=Marco Rigoli\\
\>Dipartimento di Matematica\\
\> Universita degli Studi di Milano\\
\>  Via Saldini 50\\
\> I-20133, Milano\\
\> Italy\\
\> marco.rigoli@unimi.it
\end{tabbing}}


\begin{thebibliography}{20}

\bibitem{AAR} G. Albanese, L. J. Al\'{i}as and M. Rigoli.
\emph{A general form of the weak maximum principle and some applications}.
To appear in Rev. Mat. Iberoam.

\bibitem{AHP} L. J. Al\'{i}as, A. Hurtado and V. Palmer.
\emph{Geometric analysis of Lorentzian distance function on spacelike hypersurfaces}.
Trans. Amer. Math. Soc. {\bf 362} (2010), 5083--5106.

\bibitem{AIRlorentz} L. J. Al\'{i}as, D. Impera and M. Rigoli.
\emph{Spacelike hypersurfaces of constant higher order mean curvature in
generalized Robertson-Walker spacetimes},
Math. Proc. Camb. Phil. Soc. {\bf 152} (2012), 365--383.

\bibitem{AIR} L. J. Al\'{i}as, D. Impera and M. Rigoli. \emph{Hypersurfaces of constant
higher order mean curvature in warped products}.
Trans. Amer. Math. Soc. {\bf 365} (2013), 591--621.

\bibitem{BC} J. Barbosa and A. Colares.
\emph{Stability of hypersurfaces with constant \mbox{$r$-mean} curvature},
Ann. Glob. An. Geom. {\bf 15} (1997), 277--297.

\bibitem{Ek} I. Ekeland. \emph{Nonconvex minimization problems},
Bull. Amer. Math. Soc.~{\bf 1} (1979), 443--474.

\bibitem{EGK} F. Erkekoglu, E. Garc\'\i a-R\'\i o and D. Kupeli.
\emph{On level sets of Lorentzian distance function}, Gen.\ Relativity
Gravitation {\bf 35} (2003), \mbox{1597--1615}.

\bibitem{E} M. Elbert. \emph{Constant positive $2$-mean curvature hypersurfaces}.
Illinois J. Math. {\bf 46} (2002), 247--267.

\bibitem{FS1} F. Fontenele and S. Silva. \emph{On the $m$-th mean curvature of compact
hypersurfaces}. Houston J. Math. {\bf 32} (2006), 47--57.

\bibitem{PRSbook} S. Pigola, M. Rigoli and A.G. Setti,
\textit{Vanishing and finiteness results in geometric analysis.
A generalization of the Bochner technique}. Progress in Mathematics,
{\bf 266}. Birkh\"auser Verlag, Basel, 2008.

\bibitem{Rm} A. Ranjbar-Motlagh. \emph{Rigidity of spheres in Riemannian manifolds and a
non-embedding theorem}. Bol. Soc. Brasil. Mat. {\bf 32} (2001), 159--171.

\bibitem{R} J. Roth.  \emph{Extrinsic radius pinching in space forms of nonnegative
sectional curvature}. Math. Z. {\bf 258} (2008), 227--240.

\bibitem{V} A. Veeravalli.  \emph{On the mean curvatures sharp estimates of hypersurfaces}.
Expo. Math. {\bf 20} (2002), 255--261.

\bibitem{Vl} T. Vlachos. \emph{A characterization of geodesic spheres in space forms}.
Geom. Dedicata {\bf 68} (1997), 73--78.



\end{thebibliography}
\end{document}